\documentclass[11pt]{article}

\usepackage{amsmath,makeidx,amssymb,amsthm,amscd,index,mathset,psfig,graphicx}

\newcommand{\T}{\protect\textstyle}
\newcommand{\eps}{\varepsilon}
\newcommand{\lbd}{\lambda}
\newcommand{\ipl}{\langle} 
\newcommand{\ipr}{\rangle} 

\newcommand{\Beta}{B}
\newcommand{\ga}{\tilde{\gamma}}
\newcommand{\Qt}{\widetilde{Q}}
\newcommand{\Lb}{\mathcal L}
\newcommand{\tr}{\mbox{tr}}
\newcommand{\half}{{\textstyle\frac{1}{2}}}
\newcommand{\reals}{\mathbb{R}}

\newtheorem{theorem}{Theorem}[section]

\newtheorem{lemma}[theorem]{Lemma}
\newtheorem{proposition}[theorem]{Proposition}
\newtheorem{corollary}[theorem]{Corollary}
\newtheorem{remark}[theorem]{Remark}

\begin{document}
\setcounter{footnote}{2}

\title{\bf On regularization methods for inverse problems of dynamic type}

\author{S. Kindermann$^1$ \quad  and \quad A. Leit\~ao$^{2}$ \\[2ex]
\small ${}^1$Department of Mathematics, UCLA, 520 Portalo Plaza,
Los Angeles, CA 90095  \hfill \mbox{} \\
{\small\tt \{kindermann@indmath.uni-linz.ac.at\} }  \hfill \mbox{} \\
\small ${}^2$Department of Mathematics, Federal University of
St. Catarina, P.O. Box 476,  \hfill \mbox{} \\
\small 88.040-900 Florianopolis, Brazil {\tt \{aleitao@mtm.ufsc.br\} }
\hfill \mbox{} }
\vskip0.8cm

\date{} \maketitle

\begin{abstract}
In this paper we consider new regularization methods for linear inverse
problems of dynamic type.
These methods are based on dynamic programming techniques for linear
quadratic optimal control problems.
Two different approaches are followed: a continuous and a discrete one.
We prove regularization properties and also obtain rates of convergence
for the methods derived from both approaches.
A numerical example concerning the dynamic EIT problem is used to illustrate
the theoretical results.
\end{abstract}

\section{Introduction} \label{sec:introd}

\subsubsection*{Inverse problems of dynamic type}

We begin by introducing the notion of dynamic inverse problems. Roughly
speaking, these are inverse problems in which the measuring process
--performed to obtain the data-- is time dependent. As usual, the problem
data corresponds to indirect information about an unknown parameter, which
has to be reconstructed. The desired parameter is allowed to be itself
time dependent.

Let $X$, $Y$ be Hilbert spaces. We consider the inverse problem of
finding $u: [0,T] \to X$ from the equation
\begin{equation} \label{eq:dyn-ip-c}
F(t) u(t) \ = \ y(t)\, ,\ t \in [0,T]\, ,
\end{equation}
where $y: [0,T] \to Y$ are the dynamic measured data and $F(t): X \to Y$
are linear ill-posed operators indexed by the parameter $t \in [0,T]$.
Notice that $t \in [0,T]$ corresponds to a (continuous) temporal index.
The linear operators $F(t)$ map the unknown parameter $u(t)$ to the
measurements $y(t)$ at the time point $t$ during the finite time interval
$[0,T]$. This is called a {\em dynamic inverse problem}.

If the properties of the parameter $u$ do not change during the measuring
process, the inverse problem in (\ref{eq:dyn-ip-c}) reduces to the simpler
case \ $F(t) u = y(t)$, $t \in [0,T]$, where $u(t) \equiv u \in X$. We shall
refer to this as {\em static inverse problem}.

As one would probably expect at this point, a discrete version of
(\ref{eq:dyn-ip-c}) can also be formulated. The assumption that the
measuring process is discrete in time leads to the {\em discrete
dynamic inverse problems}, which are described by the model
\begin{equation} \label{eq:dyn-ip-d}
F_k u_k \ = \ y_k\, ,\ k = 0, \dots, N
\end{equation}
and correspond to phenomena in which only a finite number of measurements
$y_k$ are available. As in the (continuous) dynamic inverse problems, the
unknown parameter can also be assumed to be constant during the measurement
process. In this case, we shall refer to this problems as {\em discrete
static inverse problems}.

Since the operators $F(t)$ are ill-posed, at each time point $t \in [0,T]$
the solution $u(t)$ does not depend on a stable way on the right hand side
$y(t)$. Therefore, regularization techniques have to be used in order to
obtain a stable solution $u(t)$. It is convenient to consider {\em time
dependent regularization techniques}, which take into account the fact
that the parameter $u(t)$ evolves continuously with the time.

In this paper we shall concentrate our attention to the (continuous and
discrete) dynamic inverse problems. The analysis of the static problems
follows in a straightforward way, since it represents a particular
subclass of the dynamic problems.

\subsubsection*{Some relevant applications}

As a first example of dynamic inverse problem, we present the {\em dynamical
source identification problem:} Let $u(x,t)$ be a solution to
$$ \Delta_x u(x,t) = f(x,t) \quad \mbox{ in } \Omega, $$
where $f(x,t)$ represents an unknown source which moves around and
might change shape with time $t$. The inverse problem in this case
is to reconstruct $f$ from single or multiple measurements of
Dirichlet and Neumann data $(u(x,t)$, $\partial_n u(x,t))$, on the
boundary $\partial\Omega$ over time $t \in [0,T]$.
Such problems arise in the field of medical imaging, e.g. brain
source reconstruction \cite{ABF02} or electrocardiography
\cite{Le03}.

Many other 'classical' inverse problems have corresponding dynamic
counterparts, e.g., the {\em dynamic impedance tomography problem}
consists in reconstructing the time-dependent diffusion coefficient
(impedance) in the equation
\begin{equation} \label{eq:dimp}
  \nabla\cdot(\sigma(.,t) \nabla ) u(.,t) \ = \ 0,
\end{equation}
from measurements of the time-dependent Dirichlet to Neumann map
$\Lambda_\sigma$ (see the review paper \cite{CIN99}). This problem
can model a moving object with different impedance inside a fluid
with uniform impedance, for instance the heart inside the body.
Notice that in this case we assume the time-scale of the movement
to be large compared to the speed of the electro-magnetic waves.
Hence, the quasi-static formulation \eqref{eq:dimp} is a valid
approximation for the physical phenomena.

Another application concerning dynamical identification problems
for the heat equation is considered in \cite{KS97,KS99}.
Other examples of dynamic inverse problems can be found in
\cite{SLDBF01, SLWV02, SVVSK01, VK89, WB95}. In particular,
for applications related to process tomography, see the conference
papers by M.H.Pham, Y.Hua, N.B.Gray; M.Rychagov, S.Tereshchenko;
I.G.Kazantsev, I.Lemahieu in \cite{Pe00}.

\subsubsection*{Inverse problems and control theory}

Our main interest in this paper is the derivation of regularization methods
for the inverse problems (\ref{eq:dyn-ip-c}) and (\ref{eq:dyn-ip-d}).
In order to obtain this regularization methods, we follow an approach based
on a solution technique for linear quadratic optimal control problems:
the so called {\em dynamic programming} which was developed in the early
50's. Among the main early contributors of this branch of optimization
theory we mention R.Bellman, S.Dreyfus and R.Kalaba (see, e.g., \cite{Bel53,
Bel57, BDS62, BK65, Dre65}).

The starting point of our approach is the definition of optimal control
problems related to (\ref{eq:dyn-ip-c}) and (\ref{eq:dyn-ip-d}).
Let's consider the following constrained optimization problem
\begin{equation} \label{eq:ccp}
  \left\{ \begin{array}{l}
      {\rm Mimimize} \ J(u,v) := \half \int_0^T \big[ \,
                     \ipl F(t)u(t)-y(t), \ L(t) [F(t)u(t)-y(t)] \ipr \\[1ex]
      \hskip4.2cm      + \ \ipl v(t), M(t) v(t) \ipr \, \big] \ dt \\[1ex]
      {\rm s.t.} \ \ u' = A(t) u + B(t) v(t), \ t \in [0,T] ,\ u(0) = u_0\, ,
   \end{array} \right.
\end{equation}
where $F(t)$, $u(t)$ and $y(t)$ are defined as in (\ref{eq:dyn-ip-c}) and
$v(t) \in X$, $t \in [0,T]$. Further, $L(t): Y \to Y$, $M(t): X \to X$,
$A(t), B(t): X \to X$ are given operators and $u_0 \in X$.
In the control problem (\ref{eq:ccp}), $u$ plays the rule of the system
trajectory, $v$ corresponds to the control variable and $u_0$ is the
initial condition. The pairs $(u,v)$ constituted by a control strategy
$v$ and a trajectory $u$ satisfying the constraint imposed by the linear
dynamic are called {\em admissible processes}.

The goal of the control problem is to find an admissible process $(u,v)$,
minimizing the quadratic objective function $J$. This is a quite well
understood problem in the literature. Notice that the objective function
in problem (\ref{eq:ccp}) is related to the Tikhonov functional for
problem (\ref{eq:dyn-ip-c}), namely
$$ \T\int_0^T \big( \| F(t) u(t) - y(t) \|_a^2
            + \alpha \| u(t) \|_b^2 \big) \ dt \, , $$
where the norms $\|\cdot\|_a$ and $\|\cdot\|_b$, as well as the
regularization parameter $\alpha > 0$, play the same rule as the
weight functions $L$ and $M$ in (\ref{eq:ccp}).

In the formulation of the control problem, we shall use as initial
condition any approximation $u_0 \in X$ for the least square solution
$u^\dag \in X$ of $F(0) u = y(0)$.
The choice of the weight functions $L$ and $M$ in (\ref{eq:ccp}) should be
such that the corresponding optimal process $(\bar u, \bar v)$ satisfies
$F(t) \bar u(t) \approx y(t)$ along the optimal trajectory $\bar u(t)$.

In order to derive a regularization method for (\ref{eq:dyn-ip-c}),
we formulate problem (\ref{eq:ccp}) for a family of operators $L_\alpha,
M_\alpha$ indexed by a scalar parameter $\alpha > 0$, and obtain the
corresponding optimal trajectories $\bar u_\alpha(t) =
\bar u_{L_\alpha,M_\alpha}(t)$. Each optimal process is obtained by using
the dynamic programming technique, where the Riccati equation (particular
case of the Hamilton-Jacobi (HJ) equation) plays the central rule. The
optimal trajectories $\bar u_\alpha(t)$ are used in order to generate a
family of regularization operators for problem (\ref{eq:dyn-ip-c}), in the
sense of \cite{EHN96}. The choice of the operators $L_\alpha$, $M_\alpha$
play the rule of the regularization parameter.
\bigskip

What concerns the discrete dynamic inverse problem (\ref{eq:dyn-ip-d}),
we define, analogous as in the continuous case, a discrete optimal control
problem of linear quadratic type
\begin{equation} \label{eq:dcp}
  \left\{ \begin{array}{l}
      {\rm Mimimize} \ J(u,v) := \sum\limits_{k=0}^{N-1}
                                 \ipl F_k u_k - y_k, L_k (F_k u_k - y_k) \ipr
                                 + \ipl v_k , M_k v_k \ipr \\
      \hspace{3.2cm}   + \ \ipl F_N u_N - y_N, L_N (F_N u_N - y_N) \ipr \\
      {\rm s.t.} \\
      u_{k+1} = A_k u_k + B_k v_k, \ k=0,\dots,N-1,\ \ u_0 \in X ,
   \end{array} \right.
\end{equation}
where $F_k$, $u_k$, $y_k$ are defined as in (\ref{eq:dyn-ip-d}) and
$v_k \in X$, $k=0,\dots,N-1$. Further the operators $L_k: Y \to Y$,
$M_k: X \to X$, $A_k, B_k: X \to X$ have the same meaning as in the
continuous optimal control problem (\ref{eq:ccp}). To simplify the
notation, we represent the processes $(u_k, v_k)_{k=1}^N$ by $(u,v)$.

Again, using the dynamic programming technique for this discrete linear
quadratic control problem, we are able to derive an iterative
regularization method for the inverse problem (\ref{eq:dyn-ip-d}). In this
discrete framework, the dynamic programming approach consists basically
of the Bellman optimality principle and the dynamic programming equation.

\subsubsection*{Literature overview and outline of the paper}

Continuous and discrete regularization methods for inverse problems have
been quite well studied in the last two decades and one can find relevant
information, e.g., in \cite{EHN96, EKN89, ES00, HNS95, Mor93, Tau94} and
in the references therein.

So far dynamic programming techniques have been mostly applied to solve
particular inverse problems. In \cite{KS97} the inverse problem of
identifying the initial condition in a semilinear parabolic equation
is considered. In \cite{KS99} the same authors consider a problem of
parameter identification for systems with distributed parameters.
In \cite{KL03}, the dynamic programming methods are used in order
to formulate an abstract functional analytical method to treat general
inverse problems.

What concerns dynamic inverse problems, regularization methods where
considered for the first time in \cite{SL02, SLWV02}. There, the authors
analyze discrete dynamic inverse problems and propose a procedure called
{\em spatio temporal regularizer} (STR), which is based on the
minimization of the functional
\begin{equation} \label{eq:louis}
\Phi(u) \ := \ \T\sum\limits_{k=0}^{N} \| F_k u_k - y_k \|_{L^2}^2 +
     \lambda^2 \T\sum\limits_{k=0}^{N} \| u_k \|_{L^2}^2 +
         \mu^2 \T\sum\limits_{k=0}^{N-1} \frac{\| u_{k+1} - u_k \|_{L^2}^2}
                                      {(t_{k+1} - t_k)^2} .
\end{equation}
Notice that the term with factor $\lambda^2$ corresponds to the classical
(spacial) Tikho\-nov-Philips regularization, while the term with factor
$\mu^2$ enforces the temporal smoothness of $u_k$.

A characteristic of this approach is the fact that the hole solution
vector $\{ u_k \}_{k=0}^N$ has to be computed at a time. Therefore,
the corresponding system of equations to evaluate $\{ u_k \}$ has
very large dimension. In the STR regularization, the associated system
matrix is decomposed and rewritten into a Sylvester matrix form. The
efficiency of this approach is based on fast solvers for the Sylvester
equation.

This paper is organized as follows:
In Section~\ref{sec:deriv} we derive the solution methods discussed in
this paper.
In Section~\ref{sec:regul} we analyze some regularization properties of
the proposed methods. 
In Section~\ref{sec:appl} we present numerical realizations of the discrete
regularization method as well as a discretization of the continuous
regularization method. For comparison purposes we consider a dynamic EIT
problem, similar to the one treated in \cite{SLWV02}.

\section{Derivation of the regularization methods} \label{sec:deriv}

We begin this section considering a particular case, namely the dynamic
inverse problems with constant operator. The analysis of this simpler
problem allow us to illustrate the dynamic programming approach followed
in this paper. In Subsections~\ref{ssec:d-dip} and \ref{ssec:d-ddip}
we consider general dynamic inverse problems and derive a continuous
and a discrete regularization method respectively.

\subsection{A tutorial approach: The constant operator case}
\label{ssec:tutor}

In this subsection we derive a family of regularization operators for the
dynamic inverse problem in (\ref{eq:dyn-ip-c}), in the particular case where
the operators $F(t)$ does not change during the measurement process, i.e.
$F(t) = F: X \to Y$, $t \in [0,T]$. The starting point of our approach is
the constrained optimization problem in (\ref{eq:ccp}). We shall consider
a very simple dynamic, which does not depend on the state $u$, but only
on the control $v$, namely: $u' = v$, $t \ge 0$. In this case, the
control $v$ can be interpreted as a {\em velocity function}. The pairs $(u,v)$
formed by a trajectory and the corresponding control function are called
{\em admissible processes} for the control problem.

Next we define the residual function $\eps(t) := F u(t) - y(t)$
associated to a given trajectory $u$. Notice that this residual function
evolves according to the dynamic
$$ \eps' = F u(t) - y'(t) = F v(t) - y'(t)\, ,\ t \ge 0\, . $$
With this notation, problem (\ref{eq:ccp}) can be rewritten in the form
\begin{equation} \label{eq:ccp-aux}
  \left\{ \begin{array}{l}
      {\rm Mimimize} \ J(\eps,v) = \half
                     \int_0^T \ipl \eps(t), L(t) \eps \ipr +
                              \ipl v(t), M(t) v(t) \ipr \ dt \\
      {\rm s.t.} \\
      \eps' = F v(t) - y'(t), \ t \ge 0 ,\ \ \eps(0) = F(0) u_0 - y(0)\, .
   \end{array} \right.
\end{equation}

The next result states a parallel between solvability of the optimal
control problem (\ref{eq:ccp}) and the auxiliary problem (\ref{eq:ccp-aux}).

\begin{proposition} \label{prop:equiv}
If $(\bar u, \bar v)$ is an optimal process for problem (\ref{eq:ccp}),
then the process $(\bar \eps, \bar v)$, with $\bar \eps :=
F \bar u(t) - y(t)$, will be an optimal process for problem
(\ref{eq:ccp-aux}). Conversely, if $(\bar \eps, \bar v)$ is an optimal
process for problem (\ref{eq:ccp-aux}), with $\eps(0) = F u_0 - y(0)$,
for some $u_0 \in X$, then the corresponding process $(\bar u, \bar v)$
is an optimal process for problem (\ref{eq:ccp}).
\end{proposition}

In the sequel, we derive the dynamic programming approach for the
optimal control problem in (\ref{eq:ccp-aux}). We start by introducing
the first Hamilton function $H: [0,T] \times X^3 \to \mathbb R$, defined by
$$  H(t,\eps,\lbd,v) \ := \ \ipl \lbd , F v \ipr - \ipl \lbd , y'(t) \ipr +
    \half [ \ipl \eps, L(t)\eps \ipr + \ipl v, M(t)v \ipr ] \, . $$
Notice that the variable $\lbd$ plays the role of a Lagrange multiplier
in the above definition. According to the Pontryagin's maximum principle,
the Hamilton function furnishes a necessary condition of optimality for
problem (\ref{eq:ccp-aux}). Furthermore, since (in this particular case)
this function is convex in the control variable, this optimality
condition also happens to be sufficient. From the maximum principle
we know that, along an optimal trajectory, the equality
\begin{equation} \label{eq:opt-cond}
0 \ = \ \frac{\partial H}{\partial v}(t,\eps(t),\lbd(t),v(t))
      \ = \ F^* \lbd(t) + M(t) v(t)
\end{equation}
holds. This means that the optimal control $\bar v$ can be obtained directly
from the Lagrange multiplier $\lbd: [0,T] \to X$, by solving the system
$$  M(t) \bar v(t) = -F^* \lbd(t)\, ,\ \forall t\, . $$
Therefore, the key task is actually the evaluation of the Lagrange
multiplier. This leads us to the HJ equation. Substituting the above
expression for $\bar v$ in (\ref{eq:opt-cond}), we can define the
second Hamilton function $\mathcal H: \mathbb R \times X^2 \to \mathbb R$
$$  \mathcal H(t,\eps,\lbd) \, := \, \min_{v \in X} \{ H(t,\eps,\lbd,v) \}
    \, = \, \half \ipl \eps , L(t) \eps \ipr - \ipl \lbd , y'(t) \ipr
    - \half \ipl \lbd , F M(t)^{-1} F^* \lbd \ipr \, . $$
Now, let $V: [0,T] \times X \to \mathbb R$ be the value function
for problem (\ref{eq:ccp-aux}), i.e.
\begin{eqnarray}
  V(t,\xi) \!\!\!\! & := & \!\!\!\!
                  \min\Big\{ \half \T\int_t^T
                  \ipl \eps(s) , L(s) \eps(s) \ipr +
                  \ipl v(s) , M(s) v(s) \ipr \, ds \ \Big| \ (\eps,v) \
                  {\rm admissible} \nonumber \\
           \!\!\! &    & \!\!\!   \label{def:cvf}
                  {\rm \ \ \ \ \ \ \ \ process\ for\ (\ref{eq:ccp-aux})\
                  with\ initial\ condition} \ \eps(t) = \xi \Big\} \, .
\end{eqnarray}
Our interest in the value function comes from the fact that this function
is related to the Lagrange multiplier $\lbd$ by: $\lbd(t) = V_\eps
(t,\bar \eps)$, where $\bar \eps$ is an optimal trajectory. From the
control theory we know that the value function is a solution of the
HJ equation
\begin{eqnarray} \label{eq:hjb}
 0 & = & V_t(t,\eps) + \mathcal H(t,\eps,V_\eps(t,\eps)) \nonumber \\
   & = & V_t + \half \ipl \eps , L(t) \eps \ipr -
         \ipl V_\eps , y'(t) \ipr -
         \half \ipl V_\eps , F M(t)^{-1} F^* V_\eps \ipr .
\end{eqnarray}
Now, making the ansatz: $V(t,\eps) = \half \ipl \eps , Q(t) \eps \ipr
+ \ipl b(t) , \eps \ipr + g(t)$, with $Q: [0,T] \to \mathbb R$, $b: [0,T] \to
X$ and $g: \mathbb R \to \mathbb R$, we are able to rewrite (\ref{eq:hjb})
in the form
\begin{multline} \label{eq:poly-2}
\half \ipl \eps , Q'(t) \eps \ipr + \ipl b'(t) , \eps \ipr + g'(t) +
\half \ipl \eps ,  L(t) \eps \ipr  -
                \ipl Q(t) \eps + b(t) ,  y'(t) \ipr \\  -  \half
\ipl Q(t) \eps + b(t) , F M(t)^{-1} F^* [Q(t)\eps+b(t)] \ipr = 0 \, .
\end{multline}
This is a polynomial equation in $\eps$, therefore the quadratic, the linear
and the constant terms must vanish. The quadratic term yields for $Q$ the
Riccati equation:
\begin{equation} \label{eq:riccati1}
Q'(t) \ = \ -L(t) + Q F M(t)^{-1} F^* Q \, .
\end{equation}
From the linear term in (\ref{eq:poly-2}) we obtain an evolution equation
for $b$
\begin{equation} \label{eq:riccati2}
b' \ = \ Q(t) F M(t)^{-1} F^* b + Q(t) y'(t)
\end{equation}
and from the constant term in (\ref{eq:poly-2}) we derive an evolution
equation for $g$
\begin{equation} \label{eq:riccati3}
g' \ = \ \half \ipl b(t), F M(t)^{-1} F^* b(t) \ipr
    + \ipl b(t) , y'(t) \ipr \, .
\end{equation}

Notice that the cost of all admissible processes for an initial
condition of the type $(T,\eps)$ is zero. Therefore we have to
consider the system equations (\ref{eq:riccati1}), (\ref{eq:riccati2}), 
(\ref{eq:riccati3}) with the final conditions
\begin{equation} \label{eq:riccati-ic}
Q(T) \ = \ 0\, , \ \ b(T) \ = \ 0\, , \ \ g(T) \ = \ 0\, .
\end{equation}
Notice that this system can be solved separately, first for $Q$, than
for $b$, and finally for $g$.

Once we have solved the initial value problem (\ref{eq:riccati1})--%
(\ref{eq:riccati-ic}), the Lagrange multiplier is given by $\lbd(t)
= Q(t) \bar\eps(t) + b(t)$ and the optimal control is obtained in the
form of the feedback control $\bar v(t) = - M^{-1}(t) F^*
[Q(t) \bar\eps(t) + b(t)]$. Therefore, the optimal trajectory of
problem (\ref{eq:ccp}) is given by
\begin{equation} \label{eq:cont-reg}
 \bar u' = - M^{-1}(t) F^* \big( Q(t)[F \bar u(t) - y(t)] + b(t) \big) \, ,\ \
 \bar u(0) = u_0 \, .
\end{equation}

By choosing appropriately a family of operators $\{ M_\alpha,L_\alpha \}_
{\alpha>0}$, it is possible to use the corresponding optimal trajectories
$\bar u_\alpha$, defined by the
initial value problem (\ref{eq:cont-reg}) in order to define a family of
reconstruction operators $R_\alpha: L^2((0,T) ; Y) \to H^1( (0,T) ; X)$, by
\begin{equation} \label{eq:cont-Ralpha}
R_\alpha(y) := u_0 - \T\int_0^t
M_\alpha^{-1}(s) F^* \big( Q(s)[F \bar u_\alpha(s) - y(s)]+b(s) \big) \ ds\, .
\end{equation}

We shall return to the operators $\{ R_\alpha \}$ in Section~\ref{sec:regul}
and prove that the family of operators defined in (\ref{eq:cont-Ralpha})
is a regularization method for (\ref{eq:dyn-ip-c}) (see, e.g., \cite{EHN96}).

\begin{remark}
It is possible to simplify the above equations to compute the optimal
trajectory $\bar u$.
If we introduce the function $\eta(t) := F^* Q(t) y(t) - F^* b(t)$, then
we can write $\bar u' = - M^{-1}(t) F^* Q(t) F \bar u + M^{-1}(t) \eta$.
Furthermore, using the equations for $Q'$ and $b'$, we have $\eta' =
-F^* L y(t) + F^* Q(t) F M^{-1}(t) \eta$. Thus, solving (\ref{eq:cont-reg})
is equivalent to solve the system
$$ \bar u' = - M^{-1}(t) F^* Q(t) F \bar u + M^{-1}(t) \eta\, ,\ \ \
   \eta'   = - F^* L y(t) + F^* Q(t) F M^{-1}(t) \eta \, . $$

This system can again be solved separately, first for $\eta$ (backwards in
time, with $\eta(T) = 0$) and then for $\bar u$ (forward in time). Notice
that the computation of both $b(t)$ and $g(t)$ is not required to build this
system. Furthermore, we do not need the derivative of the data $y(t)$.
\end{remark}

\subsection{Dynamic inverse problems}
\label{ssec:d-dip}

In the sequel we consider the dynamic inverse problem described in
(\ref{eq:dyn-ip-c}). As in the previous subsection, we shall look for
a continuous regularization strategy.

We start by considering the constrained optimization problem (\ref{eq:ccp}),
where $F(t)$, $u(t)$ and $y(t)$ are defined as in (\ref{eq:dyn-ip-c}),
$v(t) \in X$, $t \in [0,T]$, $L(t): Y \to Y$, $M(t): X \to X$,
$A(t) \equiv I: X \to X$, $B(t) \equiv 0$ and $u_0 \in X$.

Following the footsteps of the previous subsection, we define the first
Hamilton function $H: [0,T] \times X^3 \to \mathbb R$ by
$$  H(t,u,\lbd,v) \ := \ \ipl \lbd , v \ipr + \half
    [ \ipl F(t) u - y(t), L(t) (F(t) u - y(t)) \ipr + \ipl v, M(t)v \ipr ]
    \, . $$
Thus, it follows from the maximum principle: $0 = \partial H / \partial v
(t,u(t),\lbd(t),v(t)) = \lbd(t) + M(t) v(t)$, and we obtain a relation
between the optimal control and the Lagrange parameter, namely:
$\bar v(t) = - M^{-1}(t) \lbd(t)$.

As before, we define the second Hamilton function $\mathcal H: \mathbb R
\times X^2 \to \mathbb R$
$$  \mathcal H(t,u,\lbd) \, := \,
    \half \ipl F(t) u - y(t) , L(t) (F(t) u - y(t) \ipr -
    \half \ipl \lbd , M(t)^{-1} \lbd \ipr \, . $$
Since $\lbd(t) = \partial V / \partial u(t,u)$, where $V: [0,T] \times X
\to \mathbb R$ is the value function of problem (\ref{eq:ccp}), it is
enough to obtain $V$. This is done by solving the HJ equation (see
\eqref{eq:hjb})
\begin{eqnarray*}
 0 & = & V_t + \half \ipl F(t)u - y(t) , L(t) (F(t)u - y(t)) \ipr
             - \half \ipl V_u , M(t)^{-1} V_u \ipr \, .
\end{eqnarray*}
As in Subsection~\ref{ssec:tutor}, we make the ansatz $V(t,u) =
\half \ipl u , Q(t) u \ipr + \ipl b(t), u \ipr + g(t)$,
with $Q: [0,T] \to \mathbb R$, $b: [0,T] \to X$ and $g: \mathbb R \to
\mathbb R$. Then, we are able to rewrite the HJ equation above in the
form of a polynomial equation in $u$. Arguing as in (\ref{eq:poly-2}),
we conclude that the quadratic, the linear and the constant terms of this
polynomial equation must all vanish. Thus we obtain
\begin{equation} \label{eq:evol}
Q' = Q^* M(t)^{-1} Q - F^*(t) L(t) F(t) , \ \ 
b' = Q(t)^* M(t)^{-1} b + F^*(t) L(t) y(t) .
\end{equation}
%
%
The final conditions $Q(T) = 0$, $b(T) = 0$ are derived just like in the
previous subsection.%
\footnote{Since function $g$ is not needed for the computation of
the optimal trajectory, we omit the expression of the corresponding dynamic.}
Once the above system is solved, the optimal control $\bar u$ is obtained
by solving
\begin{equation} \label{eq:u}
\bar u'(t) = -M^{-1}(t) V_u(t,u) = -M^{-1}(t) [ Q(t) \bar u(t) + b(t)]
\end{equation}
with initial condition $\bar u(0) = u_0$.

Following the ideas of the previous tutorial subsection, we shall
choose a family of operators $\{ M_\alpha, L_\alpha \}_{\alpha>0}$ and
use the corresponding optimal trajectories $\bar u_\alpha$ in order to
define a family of reconstruction operators $R_\alpha: L^2((0,T) ; Y)
\to H^1( (0,T) ; X)$,
$$ R_\alpha(y) := u_0 - \T\int_0^t
                  M_\alpha^{-1}(s) [ Q(s) \bar u(s) + b(s)] \, ds \, . $$

The regularization properties of the operators $\{ R_\alpha \}$ will be
analyzed in Section~\ref{sec:regul}.

\subsection{Discrete dynamic inverse problems}
\label{ssec:d-ddip}

In this subsection we use the optimal control problem (\ref{eq:dcp}) as
starting point to derive a discrete regularization method for the
inverse problem in (\ref{eq:dyn-ip-d}).

In the framework of discrete dynamic inverse problems, we have a
trajectory, represented by the sequence $u_k$, which evolves according
to the dynamic
$$ u_{k+1} \ = \ A_k u_k \, + \, B_k v_k\, , \ k =0,1,\dots,N  $$
where the operators $A_k$ and $B_k$ still have to be chosen and
$\{ v_k \}_{k=0}^{N-1}$, is the control of the system.
As in the continuous case, we shall consider a simpler dynamic:
$u_{k+1} = u_k + v_k$, $k=0,1,\dots$ (i.e., $A_k = B_k = I$). In
the objective function $J$ of (\ref{eq:dcp}) we choose $M_k = \alpha I$,
$\alpha \in \mathbb{R}^+$, for all $k$.

In the sequel, we derive the dynamic programming approach for the
optimal control problem in (\ref{eq:dcp}). We start by introducing
the value function (or Lyapunov function) $V: \mathbb R \times X
\to \mathbb R$
$$  V(k,\xi) \ := \ \min\{ J_k(u,v)\; |\ (u,v) \in Z_k(\xi)
                    \times X^{N-k} \} \, , $$
where
\begin{multline} \label{eq:cost-funcD}
J_k(\eps,v) \ := \ \half \Big[
    \ipl F_N u_N - y_N , L_N (F_N u_N - y_N) \ipr + \\
    + \T\sum_{j=k}^{N-1} \ipl F_j u_j-y_j , L_j (F_j u_j-y_j) \ipr
    + \alpha \ipl v_j , v_j \ipr \Big]
\end{multline}
and $Z_k(\xi) := \{ u \in X^{N-k+1}\; |\ u_k = \xi ,\ u_{j+1} =
u_j + v_j,\ j=k,\dots,N-1 \}$.
(Compare with the definition in (\ref{def:cvf})). The Bellman
principle for this discrete problem reads
\begin{equation} \label{eq:hjb-d}
V(k,\xi) \ = \ \min_{v \in X} \{ V(k+1, \xi + v)
+ \half \ipl F_k \xi - y_k , L_k (F_k \xi - y_k) \ipr +
  \T\frac{\alpha}{2} \ipl v , v \ipr \} \, .
\end{equation}
The optimality equation (\ref{eq:hjb-d}) is the discrete counterpart of
the HJ equation (\ref{eq:hjb}). Now we make the ansatz for
the value function: $V(k,\xi) = \half \ipl \xi , Q_k \xi \ipr +
\ipl b_k , \xi \ipr + g_k$.
Notice that the value function satisfies the boundary condition:
$V(N,\xi) = \half \ipl F_N \xi - y_N, L_N (F_N \xi - y_N) \ipr$.
Therefore,
\begin{equation} \label{eq:discbc}
Q_N = F_N^* L_N F_N  \quad  b_N = - F_N^* L_N y_N .
\end{equation}

As in the continuous case, the optimality equation has to be solved
backwards in time ($k = N-1, \dots, 0$) recursively. A straightforward
calculation shows that the minimizer of (\ref{eq:hjb-d}) is given by
$\bar v = -(Q_{k+1} + \alpha I)^{-1} (Q_{k+1} \xi + b_{k+1})$. Substituting
in (\ref{eq:hjb-d}), we obtain a recursive formula to compute $Q_k$,
$b_k$ and $u_k$:
\begin{eqnarray}
Q_{k-1} & = & \alpha (Q_k + \alpha I)^{-1} Q_k + F_{k-1}^* L_{k-1} F_{k-1}
              \quad k = N \ldots 2 \label{q} \\
b_{k-1} & = & \alpha (Q_k + \alpha I)^{-1} b_k - F_{k-1}^* L_{k-1} y_{k-1}
              \quad k = N \ldots 2 \label{b} \\
u_{k+1} & = & (Q_{k+1} + \alpha I)^{-1} (\alpha u_k - b_{k+1})
              \quad k = 0 \ldots N-1 \label{u}
\end{eqnarray}
Together with the end conditions (\ref{eq:discbc}) and an arbitrary
initial condition $u_0$, these recursions can be solved backwards for
$Q_k$, $b_k$ and forwards for $u_k$.
In the sequel we verify that the iteration in \eqref{q}, \eqref{b},
\eqref{u} is well defined.

\begin{lemma}
The recursion \eqref{q} with the condition \eqref{eq:discbc},
defines a sequence of self-adjoint positive semi-definite operators $Q_k$.
In particular, $(Q_k + \alpha I)^{-1}$ exists and is bounded for all $k$.
Moreover,
$$ \|Q_k\| \leq \alpha + \max_k \|F_k\|^2. $$
\end{lemma}
\noindent {\it Proof:}
Since a sum of two bounded selfadjoint operators is symmetric, it
follows by induction that $Q_k$ are self-adjoint for all $k$.
Denote by $\sigma(Q_k)$ its spectrum, we can prove by induction that
$$ \sigma(Q_k) \subset
\left[ 0, \alpha + \max_k \|F_k\|^2 \right]. $$
Indeed, if $Q_{k+1}$ has this property, then $(Q_{k+1} + \alpha I)^{-1}$
exists, and
$$ B_{k+1} :=\left( \alpha (Q_{k+1} + \alpha I)^{-1} Q_{k+1} \right) $$
is positive semidefinite and bounded by $\|B_{k+1}\| \leq \alpha$. Hence,
by the minimax characterization of the spectrum we obtain
\begin{eqnarray*}
\sigma(Q_k) & \geq & \inf_{\|x\|\leq 1} (x,Q_k x)
              \geq   \inf_{\|x\|\leq 1} (x,B_k x) + \inf_{\|x\|\leq 1}
                     (x,F^*_k F^*_k x) \geq 0 \\
\sigma(Q_k) & \leq & \sup_{\|x\|\leq 1} (x,Q_k x)
              \leq   \sup_{\|x\|\leq 1} (x,B_k x) + \sup_{\|x\|\leq 1}
                     (x,F^*_k F^*_k x)
              \leq   \alpha + \|F_k\|^2,
\end{eqnarray*}
concluding the proof. \hfill $\Box$

\section{Regularization properties} \label{sec:regul}

Before we examine the regularization properties of the methods derived in
Section~\ref{sec:deriv}, let us state a result about existence and uniqueness
of the Riccati equations~\eqref{eq:evol}.

\begin{theorem}
If $F$, $L$, $M \in C([0,T],\Lb(X,Y))$, then the Riccati equation
\eqref{eq:evol} has a unique symmetric positive semidefinite solution
in $C^1([0,T],\Lb(X))$.
%
\end{theorem}

\noindent {\it Proof:} In \cite{Ba81} the uniqueness and 
positivity of a weak solution to \eqref{eq:evol} in the form

\begin{equation} \label{eq:weak}
\Qt(t) = \T\int_t^T \Qt(s)^* M^{-1}(s) \Qt(s) - F(s)^* L(s) F(s) ds,
\end{equation}
is proven. If $F, L, M$ are continuous then, by Lebesgues 
Theorem, $\Qt$ is continuously differentiable, and hence a strong 
solution. The symmetry of $\Qt$ follows from the uniqueness, 
since $\Qt^*$ satisfies the same equation as $\Qt$. Existence of 
a solution to \eqref{eq:evol}, \eqref{eq:u} is standard, as these 
are linear equations (cf. \cite{Sh97}). \hfill$\Box$
%

\begin{remark}
It is well known in control theory that the existence of a solution to
\eqref{eq:evol} can be constructed from the functional
\begin{multline} \label{eq:value}
V(t,\xi) := 
            \min_{\substack{ {u(t) = \xi} \\ {u \in H^1([t,T],X)} } }
            \half
            \T\int_t^T \ipl F(s) u(s) - y(s), L(s)[F(s) u(s) - y(s)] \ipr
            \\ + \ipl u'(s), M(s) u'(s) \ipr ds.
\end{multline}
This functional is quadratic in $u$ and, from the Tikhonov regularization
theory (see, e.g., \cite{EHN96}), it admits a unique solution $u$, and is
quadratic in $\xi$. Furthermore, the leading quadratic part $(\xi , Q(t)\xi)$
is a solution to the Riccati Equation.
\end{remark}

Next we consider regularization properties of the method derived in
Subsection~\ref{ssec:d-dip}. The following lemma shows that the solution
$u$ of \eqref{eq:u} satisfies the necessary optimality condition for the
functional
\begin{equation} \label{eq:cost-func}
 J(u) \ = \ \half
            \T\int_0^T \ipl F(s) u(s) - y(s), L(s)[F(s) u(s) - y(s)] \ipr
            \\ + \ipl u'(s), M(s) u'(s) \ipr ds
\end{equation}
(notice that this is the cost functional $J(u,v)$ in \eqref{eq:ccp} with
$v = u'$).

\begin{lemma}
Let $Q(t)$, $b(t)$, $u(t)$ be defined by \eqref{eq:evol}, \eqref{eq:u},
together with the boundary conditions $Q(T) = 0$, $b(T) = 0$ and
$u(0) = u_0$. Then, $u(t)$ solves
\begin{equation} \label{eq:diff}
F^*(t) L(t) F(t) u(t) - M(t) u(t)'' = F^*(t) L(t) y(t),
\end{equation}
together with the boundary conditions $u(0) = u_0$, $u'(T) = 0$.
\end{lemma}
\noindent {\it Proof:}
Equation \eqref{eq:diff} follows form equations
\eqref{eq:evol}, \eqref{eq:u} by differentiation:
\begin{eqnarray*}
-M(t) u''(t) & = & \tfrac{d}{d t} (Q u(t) + b(t)) = Q(t)' u(t)
                   + b'(t) + Q(t) u'(t) \\
             & = & -F(t)^* L(t) F(t) u(t) + F(t)^* L(t) y(t)
 \end{eqnarray*}
The boundary condition $u(0) = u_0$ holds by definition and the
identity $u'(T) = 0$ follows from \eqref{eq:u} and the boundary
conditions for $Q$ and $b$. \hfill$\Box$
\bigskip

Since the cost functional in \eqref{eq:cost-func} is quadratic, the
necessary first order conditions are also sufficient. Thus, the solution
$u(t)$ of \eqref{eq:diff} is actually a minimizer of this functional.
Including the boundary conditions we obtain the following corollary:

\begin{corollary} \label{cor:cont-min}
The solution $u(t)$ of \eqref{eq:diff} is a minimizer of the Tikhonov
functional in \eqref{eq:cost-func} over the linear manifold
$$ \mathcal H:= \{ u \in H^1([0,T],X) \ | \ u(0) = u_0 \} \, . $$
\end{corollary}

In particular, this means that the above procedure is a regularization
method for the inverse problem \eqref{eq:dyn-ip-c}.
Bellow we summarize a stability and convergence result. The proof uses
classical techniques from the analysis of Tikhonov type regularization
methods (cf. \cite{EHN96}, \cite{EKN89}) and thus is omitted.

\begin{theorem}
Let $M(t) \equiv \alpha I$, $\alpha > 0$, $L(t) > 0$, $t \in [0,T]$
and $J_\alpha$ be the corresponding Tikhonov functional given by
\eqref{eq:cost-func}. \\
\underline{Stability:} Let the data $y(t)$ be noise free and denote by
$u_\alpha(t)$ the minimizer of $J_\alpha$.
Then, for every sequence $\{\alpha_k\}_{k \in \mathbb{N}}$ converging to zero,
there exists a subsequence $\{\alpha_{k_j}\}_{j \in \mathbb{N}}$, such that
$\{ u_{\alpha_{k_j}} \}_{j \in \mathbb{N}}$ is strongly convergent. Moreover,
the limit is a minimal norm solution. \\
\underline{Convergence:} Let $\|y^\delta(t) - y(t)\| \leq \delta$.
If $\alpha = \alpha(\delta)$ satisfies
$$ \lim_{\delta \to 0} \alpha (\delta) = 0 \mbox{ \ \ and \ \ }
   \lim_{\delta \to 0} \delta^2 / \alpha (\delta) = 0 \, . $$
Then, for a sequence $\{\delta_k\}_{k\in\mathbb{N}}$ converging to zero, there
exists a sequence $\{\alpha_k := \alpha(\delta_k)\}_{k \in \mathbb{N}}$
such that $u_{\alpha_k}$ converges to a minimal norm solution.
\end{theorem}

A result similar to the one stated in Corollary~\ref{cor:cont-min}
holds for the discrete case:

\begin{lemma} Let $Q_k,b_k,u_k$ be defined by \eqref{q}, \eqref{b}
and \eqref{u}, together with the boundary conditions \eqref{eq:discbc}.
Then $u_k$ satisfies
\begin{equation}\label{eq:diffdisc}
F_k^* L_k F_k u_k - \alpha \left(u_{k-1} - 2 u_k +
u_{k+1} \right) = F_k^* L_k y_k, \quad k= 1\ldots n
\end{equation}
together with the boundary condition $u(0) = u_0$, $u_{n+1} = u_n$.
\end{lemma}

Equation~\eqref{eq:diffdisc} is the necessary (and by convexity also
sufficient) condition for a minimizer of $J_0$ in \eqref{eq:cost-funcD}.
This proves the following corollary:

\begin{corollary}
The sequence $u_k$ is a minimizer of the Tikhonov functional
\eqref{eq:cost-funcD} over all $(w_k)$ with $w_0 = u_0$.
\end{corollary}

\section{Application to dynamic EIT problem} \label{sec:appl}

After a spacial discretization of the operator equation \eqref{eq:dyn-ip-c},
the differential equations \eqref{eq:evol}, \eqref{eq:u} can be solved by
standard methods for ordinary differential equations, such as the
Euler-Method or Runge-Kutta-Methods.
Choosing $M(t) \equiv Id:X \to X$ and $L(t) \equiv \alpha^{-1} Id:Y \to Y$
in \eqref{eq:evol}, \eqref{eq:u} we obtain
\begin{eqnarray*}
 Q'(t) &=& - \alpha^{-1} F(t)^* F(t) +  Q(t)^* Q(t) \\
 b'(t) &=&  Q(t)^* b(t)  + \alpha^{-1} F(t)^* y(t) \\
 u'(t) &=& - Q u(t) - b(t)
 \end{eqnarray*}
From a computational point of view, the first of these is the most
expensive one, as it is nonlinear and involves matrix products.
Once $Q(t)$ is known, the equations for $b,u$ are linear and only
involve matrix-vector multiplications.

The simplest approach is to use an explicit Euler method for solving
the equation for $Q$ backwards in time ($t_k = \tfrac{k}{n_T} T$,
$\Delta t = \tfrac{1}{n_T} T$).
\begin{eqnarray} \label{explicit}
Q_{k-1} \!\!\!\!\!&=& \!\!\!\!\!
          Q_k - \Delta t \left( - \alpha^{-1} F(t_k)^*F(t_k) +
          Q_k(t)^*Q_k(t)\right)  \  k = n-1,\ldots, 0 \\
b_{k-1} \!\!\!\!\!&=& \!\!\!\!\!
          b_k - \Delta t \left( Q_k b_k + \alpha^{-1} F(t_k)^*
          y(t_k)\right)  \\
u_{k+1} &=& u_k + \Delta t \left( - Q_k u_k - b_k \right) ,
\end{eqnarray}
with $Q_n = 0$, $b_n = 0$.
It is well known, that an explicit method is conditionally stable.
The iteration matrix for \eqref{explicit} is $ (I - \Delta t Q_k)$.
An analogy to Landweber iteration \cite{EHN96} a stability criterion
is that
\begin{equation} \label{stable}
\Delta t \leq \| Q_k \|^{-1} .
\end{equation}
This condition is satisfied if $\Delta t$ is small enough, as the
following Theorem states:

\begin{theorem}
Let the following CFL-condition be satisfied
\begin{equation}\label{eq:CFL}
\alpha^{-1} (\Delta t)^2 \max_{t\in [0,T]} \|F(t)\|^2 \leq \half .
\end{equation}
Then \eqref{explicit} defines a sequence of positive definite selfadjoint
operators $Q_k$ such that \eqref{stable} hold.
\end{theorem}
\noindent {\it Proof:}
It is trivial that $Q_{k-1}$ is selfadjoint if $Q_k$ is.
The iteration can be written as
$$ \Delta t Q_{k-1} = ( I - \Delta t Q_k) \Delta t Q_k +
   \alpha^{-1} (\Delta t)^2 F_k^*F_k. $$
If the spectrum $\sigma$ of $Q_k$ satisfies $\sigma(\Delta t Q_k)
\subset [0,1]$, then the right hand side of the iteration is a sum
of two positive definite operators and hence the left hand side is
also positive definite. Moreover,
$$ \|\Delta t Q_{k-1}\| \leq
   \half + \alpha^{-1} (\Delta t)^2 \|F_k\|^2 . $$
If $\alpha^{-1} (\Delta t)^2 \|F_k\|^2 \leq \half$ holds, then we obtain
by induction that $\sigma(\Delta t Q_{k-1}) \subset [0,1]$ for all $k$,
which implies \eqref{stable}. \hfill$\Box$
\bigskip

If follows from the last theorem that $\Delta t$ has to be chosen
proportional to $\sqrt{\alpha}$. If the regularization parameter
is small, this requires very  small time-steps. In this case an
alternative is to use the discrete versions
\eqref{q},\eqref{b},\eqref{u}, which are quite similar to an
implicit Euler schema. Contrary to the explicit Euler steps, it
does not require any restriction on $\Delta t$.

In this section, we apply our regularization method to a dynamic
inverse problem, namely the linearized impedance tomography problem,
i.e. one is faced with the problem of determining a time-dependent
diffusion coefficient $\ga(x,t)$ in the equation
\begin{equation} \label{pde}
 \nabla.\left(\ga(.,t) \nabla u\right) = 0 \quad \mbox{ in }  \Omega
 \end{equation}
from the Neumann-to-Dirichlet operator: 
$$ \Lambda_{\ga}: \tfrac{\partial}{\partial n} u |_{\partial\Omega} \to
   u |_{\partial\Omega}  \qquad
   u \mbox{ solution to the Neumann problem \eqref{pde}} . $$
We consider $\Lambda_{\ga}$ an operator mapping a subspace $L^2(\partial
\Omega)$ into itself. Since the Neumann data have to satisfy the compatibility
condition $\int_{\partial \Omega} \tfrac{\partial}{\partial n} u = 0$, the
domain of definition of $\Lambda_{\ga}$ has to incorporate this condition.
It is well known (see, e.g., \cite{Is98}) that $\Lambda_{\ga}$ is a compact
operator between Hilbert-spaces, hence we can consider it an element of the
space of Hilbert-Schmidt operators $H$ and use the Hilbert-Schmidt norm on
this space. The parameter-to-data operator can be written as
$F: X \subset L^2([0,T]\times\Omega) \to H$, $F(\ga) := \Lambda_{\ga}$.

The subset $X$ is the set of $\ga$ such that $\ga$ is bounded from below
and above by positive constants, which is necessary to ensure ellipticity
of \eqref{pde}.
Since the operator $F$ is nonlinear, for a successful application of the
dynamic algorithm we will consider a linearization around $1$, using
$F(\ga) - F(1) \sim F'(1) (\ga - 1)$. Notice that $F(1)$ can be computed
{\em a priori}, therefore we consider the data to our problem to be
$y = F(\ga) - F(1)$ and the corresponding unknown $\gamma(x,t) =
\ga(x,t) - 1$. This gives the linearized problem
$$ F'(1) \gamma = y , $$
where $\gamma$, $y$ both depend on time. Hence, we can solve this problem
within the framework developed in Subsection~\ref{ssec:tutor}.

\subsection{Discretization}
We briefly comment about the discretization of the Neumann-to-Dirichlet
operator. We use piecewise linear finite element functions on the boundary:
\ $X_b:=\{ \sum_i g_i \tilde{\phi}_i(x) | \ x \in \partial\Omega\}$.
The functions $\phi_i$ are the boundary-trace of the well-known
Courant-element functions. Equation~\eqref{pde} is also solved by finite
elements. Let $\phi_i$ be the piecewise linear and continuous ansatz
functions on a triangular mesh. These ansatz functions form the basis
for the finite-element space to solve \eqref{pde} and also for the
discretization of the space $X$, i.e. $\gamma$ is represented in the
discrete setting by a sum of $\phi_i$. If the Neumann data are in $X_b$,
i.e. $\tfrac{\partial}{\partial n} u = \sum_i g_i \tilde{\phi_i}$, then
equation \eqref{pde} corresponds to a discrete linear equation of the form
$$ \left( \begin{array}{cc} A_{11} & A_{12} \\
 A_{21} & A_{22} \end{array} \right) \left(\begin{array}{c}
  u_i \\ u_b \end{array} \right)  =
  \left(\begin{array}{c}
  0 \\ M g \end{array}\right), $$
where the matrices $A_{11},A_{12},A_{21},A_{22}$ are sub-matrices of the
stiffness matrix $A_{i,j} = \int_\Omega \gamma \nabla \phi_i \nabla \phi_j$
with  respect of a splitting of the indices into the interior and boundary
components. The matrix $M$ is coming from the contribution of the
Neumann-data in the discretized  equations:
\begin{equation} \label{bound}
 M_{i,l} = \T\int_{\partial \Omega} \tilde{\phi_l} \tilde{\phi}_i d\sigma  
\end{equation}
In order to deal with the compatibility condition we specify a reference
boundary index $i_*$ and set $g_{i^*}= 0$. The corresponding rows and
columns in the matrices are canceled out. The variables connected with
interior points can be eliminated from the discrete equation by taking
the Schur-Complement, this gives the matrix
\begin{equation} \label{Gdef}
 G := \left( A_{22} - A_{21} A_{11}^{-1} A_{12} \right)^{-1} M. 
\end{equation}
This matrix corresponds to a mapping $\tilde{\Lambda}:X_b^* \to  X_b^*$, with
$$ X_b^*:= \{ \T\sum_i g_i \tilde{\phi}_i(x) | \quad g_{i^*} = 0,
   \quad x \in \partial\Omega\} . $$
Identifying the space $X_b$ with the $\reals^n$ via
$\sum_i g_i \tilde{\phi}_i(x) \Leftrightarrow (g_i)$, the discrete
Neumann-to-Dirichlet operator is represented on $\reals^n$ by
multiplication of the matrix $G$.

We calculate the Hilbert-Schmidt inner product for discrete
Neumann-to-Dirichlet operators $\Lambda_1$, $\Lambda_2$ coming
from the above discretizations. These operators have the form
$\Lambda_k \tilde{\phi}_i \to \T\sum_l ({G_k})_{l,i} \tilde{\phi}_l$,
$k \in \{1,2\}$, where $G_k$ is as in \eqref{Gdef}, corresponding to
different coefficients $\gamma$. Note that $G_k$ can be written as
$G_k = S_k M$, where $S_k$ is a symmetric matrix and $M$ the boundary
mass matrix \eqref{bound}.

The Hilbert Schmidt inner product is defined as
$(\Lambda_1,\Lambda_2) = \sum_i (\Lambda_1 e_i, \Lambda_2 e_i)$,
where $e_i$ is a orthonormal basis and $(.,.)$ is the usual $L^2$ inner
product. In our case we chose $(e_i)$ orthonormal such that $\mbox{span}(e_i)
= \mbox{span} (\tilde{\phi}_i)$. Each basis can be transformed into each
other:
$\tilde{\phi}_i = \T\sum_k \beta_{i,k} e_k$,
$e_k = \sum_l \gamma_{k,l} \tilde{\phi}_{l}$.

Denote by $B$, $\Gamma$ the matrices: $B = (\beta_{i,k})$, $\Gamma :=
(\gamma_{k,l})$. From the orthogonality of $(e_i)$ the following identities
can be derived: $M = B B^T$, $\Gamma = \Beta^{-1}$.
Now $\Lambda e_k$ is given by $\Lambda e_k = \sum_{l} A_{k,l}
\tilde{\phi}_l$ and further $A = \Gamma G^T$.

Finally the Hilbert-Schmidt inner product can be calculated to
($\tr$ denotes the trace of a matrix):
$$ (\Lambda_1,\Lambda_2) = \tr(\Gamma G_1^T    M ( \Gamma G_2^T)^T) =  
\tr(G_1^T M G_2 \Gamma^T \Gamma ) =  \tr(M S_1^T M S_2 M M^{-1} )  $$
$$ = \tr(M S_1^T M S_2 )  =  \tr( S_1^T M S_2 M) = \tr(G_1 G_2), $$
where we used $ \Gamma^T \Gamma = M^{-1}$, and the symmetry of $S,M$, and the
identity $\tr(A B) = \tr(B A)$.

\subsection{Numerical Results}

\begin{figure}[t]
  \begin{center}
  \includegraphics[width=12cm]{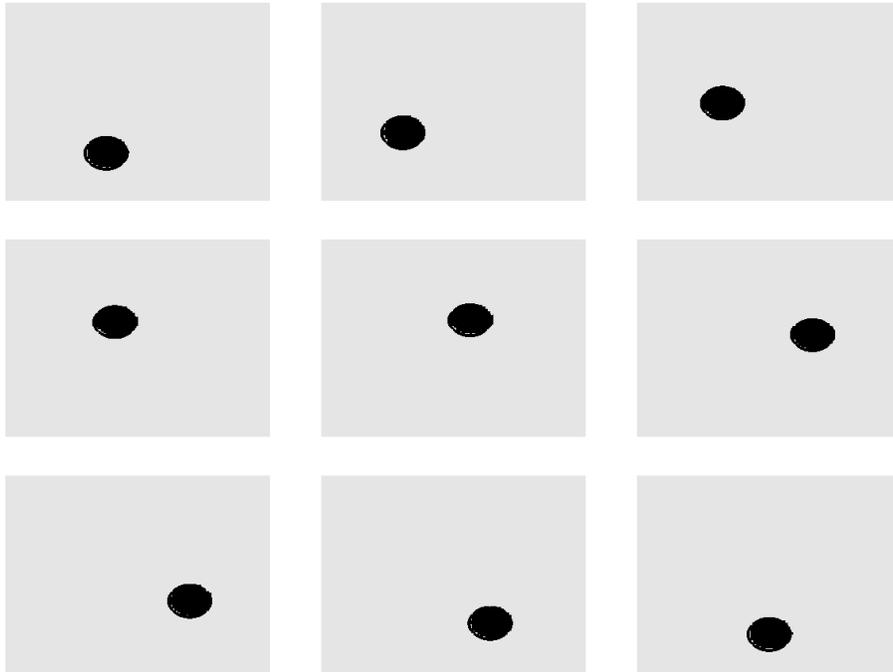} \label{fig:exact}
  \caption{Exact solution of the dynamic EIT problem.}
  \end{center}
\end{figure}

As test example for the linearized impedance tomography problem
we considered equation \eqref{pde} on a unit square: $\Omega = [0,1]^2$.
As conductivity $\gamma(x,t)$ we used a piecewise constant function,
with support on a moving circle:
$$ \gamma(x,t) := 1 + 2 \chi_{B_{x_t,0.08}}, $$
here $\chi$ denotes the characteristic function, $B_{x,r}$ denotes
a circle with center at $x$ and radius $r$. The time-varying center
is chosen as
$$ x_t := \left( \begin{array}{l}0.4 - 0.2 \cos( 2 \pi t) \\
   0.5 - 0.2\sin(2 \pi t) \end{array}\right).\qquad t \in [0,1] $$
and is shown in Figure\ref{fig:exact}.

For the computation we used a uniform discretization, with 25 subdivisions
of the interval $[0,1]$ in each coordinate direction. The data are sampled
at $t_i = \tfrac{i}{50}$ using $51$ uniform distributed sample points of
the interval $[0,1]$.

We experimented both with the explicit Euler algorithm and the discrete
version. However the first one has the drawback of needing a CFL condition
\eqref{eq:CFL}. For small $\alpha$ this requires a very fine discretization
of the time-interval, which makes the method not very feasible. Hence for
the numerical results we used the discrete version, which is free of a CFL
condition. 

\begin{figure}[ht]
  \begin{center}
  \includegraphics[width=12cm]{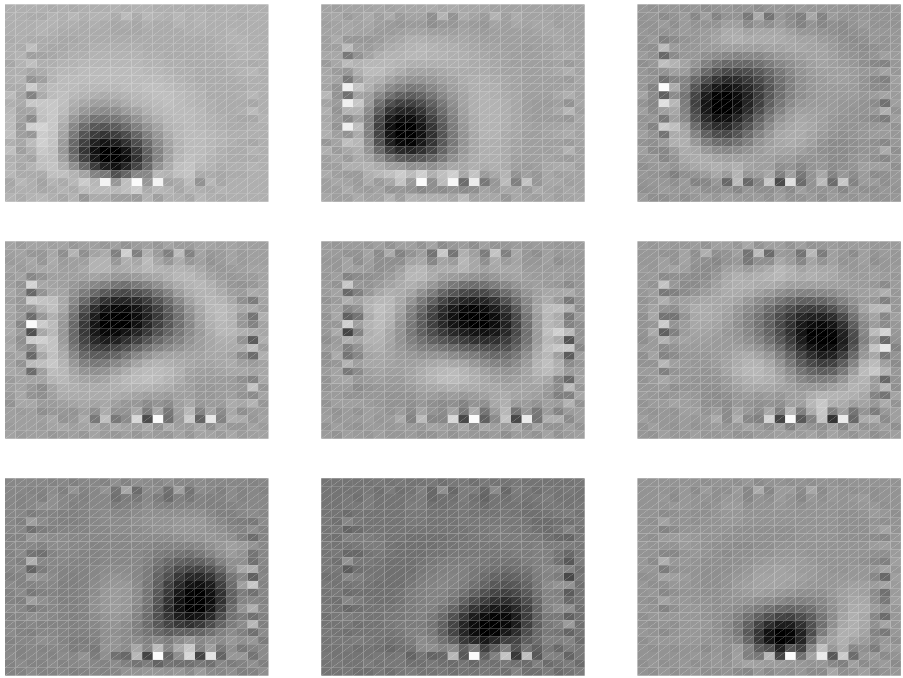} \label{fig:resprime}
  \caption{Reconstruction result for linearized data without noise.}
  \end{center}
\end{figure}

For the first example we simulated data for the linearized problem, i.e.
$$ y = F'(1) (\gamma -1). $$
The data were computed on a finer unstructured grid, in order to avoid
inverse crimes. In Figure~\ref{fig:resprime} we show a density plot of
the results for different time-points.

\begin{figure}[ht]
  \begin{center}
  \includegraphics[width=12cm]{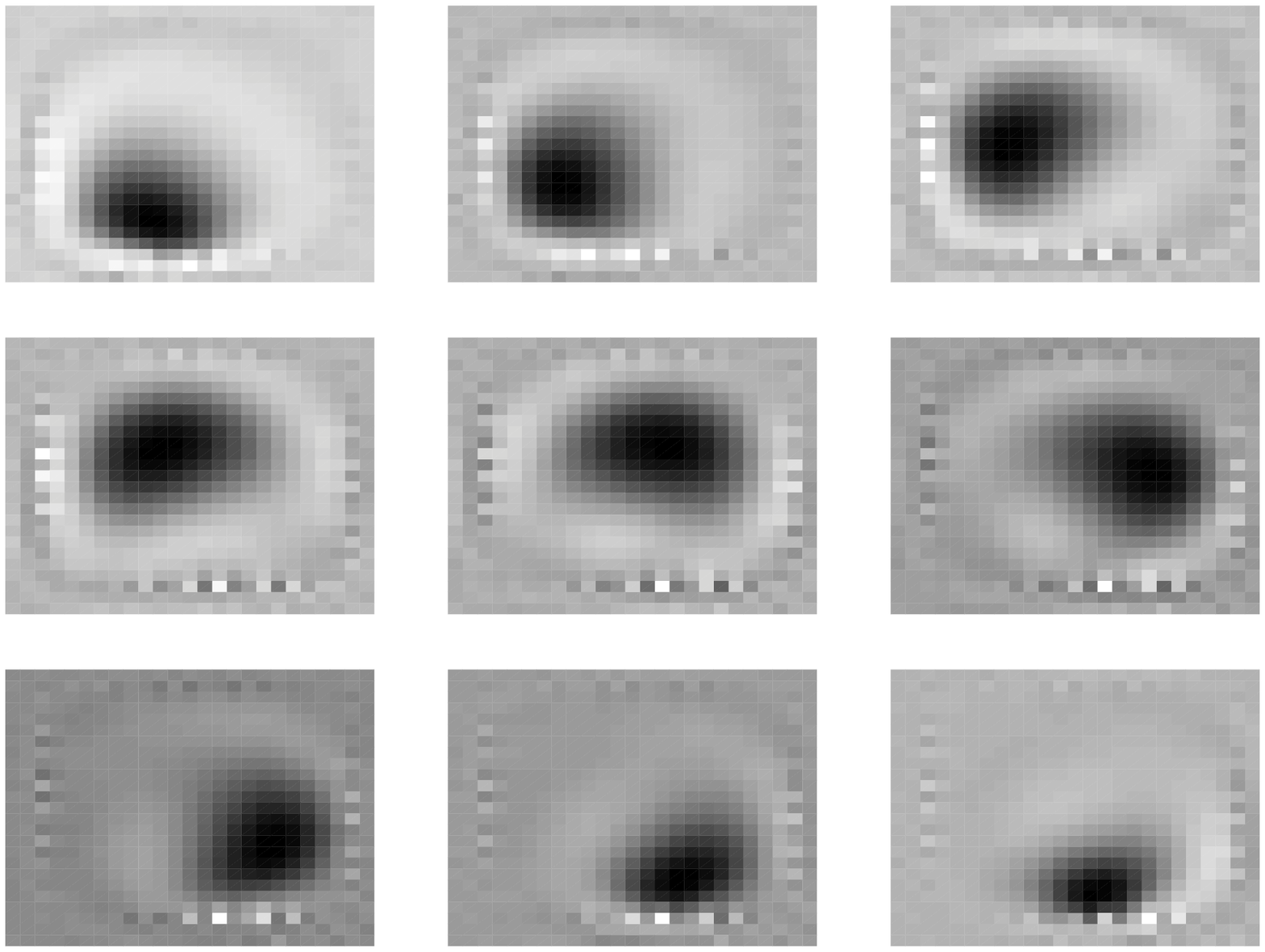} \label{fig:resnonlin}
    \caption{Reconstruction result for nonlinear data with 5\% random noise.}
    \end{center}
\end{figure}

For the second example we used nonlinear data
$$ y = F(\gamma) - F(1). $$
Again we computed this on a finer grid. Additionally, we added 5\% random
noise. Thus, we have in this case both an error due to noise and a
systematic error coming from the fact that we used a linearized model
for data corresponding to a nonlinear problem.
Figure~\ref{fig:resnonlin} shows the result for this case.

\section{Conclusions}

Each method derived in this paper require, in a first step, the solution
of an evolutionary equation (of Hamilton-Jacobi type). In a second step,
the components of the solution vector $\{ u_k \}$ are computed one at a
time. This strategy reduces significantly both the size of the systems
involved in the solution method, as well as storage requirements needed
for the numerical implementation. These points turn out to become critical
for long time measurement processes.

Some detailed considerations about complexity:
Assume that all $F(t_k)$ are discretized as $(n\times m)$ matrices.
The main effort is the matrix multiplication for the update step for
$Q_k$: In each step this requires $\mathcal{O}(n^3 + n^2 m)$ calculations.
Hence the overall complexity is of the order $\mathcal{O}(n_T(n^3 + n^2m))$
operations. If the discrete version is used, then in each step a
matrix-inversion has to be performed, which is also of the order
$\mathcal{O}(n^3)$. which leads to the same complexity as above.
In contrast, the method in \cite{SL02} requires
$\mathcal{O}((n+n_T)^3 + (n_T+m) n n_T) $. Although this is only
of cubic order in comparison to a quartic order complexity for the
dynamic programming approach, it is cubic in $n_T$. Hence if $n_T$
is large, the method proposed in this paper (which is linear in
$n_T$) will be more effective than the method in \cite{SL02}.

The numerical results show the feasibility and the stability of our
method. Note that the results are more smeared out at the center of the
square, which is clear since the identification problem is less stable
if the boundary is further away.

\section*{Acknowledgments}
The work of S.K. is supported by Austrian Science Foundation under
grant SFB F013/F1317 and by NSF grant Nr. DMI-0327077.
S.K. is on leave form the Industrial Mathematics 
Institute, Johannes Kepler University Linz, Austria.

Part of this paper was written during a sabbatical stay of A.L. at
RICAM Institute (Linz).
A.L. acknowledge support of the Austrian Academy of Sciences and of
CNPq under grants 305823/2003-5 and 478099/2004-5.

The authors would like to thank Prof K. Kunisch (Graz) for the fruitful
discussions about optimal control and optimization theory.


\end{document}